\documentclass[12pt]{amsart}
\newfont{\sheaf}{eusm10 scaled\magstep1}

\newcommand{\C}{\ensuremath{\mathbb{C}}}

\newcommand{\Z}{\ensuremath{\mathbb{Z}}}
\newcommand{\F}{\ensuremath{\mathbb{F}}}
\newcommand{\Q}{\ensuremath{\mathbb{Q}}}

\newcommand{\hol}{\ensuremath{\mathcal{O}}}

\newcommand{\PP}{\ensuremath{\mathbb{P}}}

\newcommand{\Proof}{{\it Proof. }}
\newcommand{\ra}{\ensuremath{\rightarrow}}

\newcommand{\sW}{{\mathcal W}}
\newcommand{\sS}{{\mathcal S}}

\newcommand{\sG}{{\mathcal G}}

\newcommand{\sK}{{\mathcal K}}

\newcommand{\la}{\lambda}

\newcommand{\Ga}{\Gamma}
\newcommand{\De}{\Delta}
\newcommand{\de}{\delta}

\newcommand{\ga}{\gamma}

   \newtheorem{defn}{Definition}
\newtheorem{thm}{Theorem}

\newtheorem{prop}[thm]{Proposition}
\newtheorem{rem}{Remark}
\newtheorem{ex}{Example}
\newtheorem{lem}[thm]{Lemma}
\newtheorem{cor}[thm]{Corollary}

\def\Bbb{\bf}
\def\P{{\Bbb P}}

\def\eea{\end{eqnarray*}}
\def\bea{\begin{eqnarray*}}

\title{SURFACE CLASSIFICATION AND LOCAL AND GLOBAL
FUNDAMENTAL 	GROUPS, I .}

\author{Fabrizio Catanese
   Universit\"at  Bayreuth\\}

\begin{document}
{This article is dedicated to 
 Guido Zappa, the sweet (grand-?)
father of Italian Algebra
   and Geometry, on occasion of his 90-th birthday\footnote{
The research of the  author was performed in the realm  of the
EAGER EEC, and of the Schwerpunkt 'Globale Methoden in der
komplexen Geometrie'.}.}

\begin{abstract}
Given a smooth complex surface $S$, and a compact connected global normal
crossings divisor $ D = \cup_i D_i$, we consider the local
fundamental group $\pi_1 (T \setminus D) $, where $T$ is
a good tubular neighbourhood of $D$.

One has
an exact sequence  $1 \ra \sK \ra \Ga : = \pi_1 (T - D) \rightarrow 
\Pi : = \pi_1
(D)  \ra 1$, and the kernel $\sK$ is normally generated by
geometric loops $\ga_i$ around the curve $D_i$. Among the main results,
which are strong generalizations of a  well known theorem of Mumford, is
the nontriviality of $\ga_i$ in  $\Ga = \pi_1 (T - D)$, provided all
the curves $D_i$ of genus zero have selfintersection $D_i^2 \leq -2$
(in particular this holds if
the canonical divisor $K_S$ is nef on $D$), and under the technical assumption
that the dual graph of $D$ is a tree.

    \end{abstract}

\date{January 5, 2006}
\maketitle


\section{Introduction}

In his first mathematical paper \cite{mu61} David Mumford solved the
conjecture of Abhyankar showing that, over the complex numbers
$\C$, a normal singular point $P$ of an algebraic surface $X$
   is indeed
   a smooth point if and only if it is topologically simple :
   more precisely, if and only if the local fundamental group
$\pi_{1  , \ loc} (X  , P )$ is trivial.

He derived from this result  the interesting Corollary
that the local ring $\hol_{X,P}$ of a normal singular point
is factorial if and only if either $P$ is a smooth point,
or $\pi_{1  , \ loc}(X  , P )$ is the binary icosahedral group,
and the singularity is then analytically isomorphic
to
$$ \{ (x,y,z) \in \C^3 | z^2 + x^3 + y ^5 = 0 \}$$
(a shorter  independent proof  of this corollary was later
found by Shepherd-Barron, cf. \cite{s-b99}:  this proof is 
  similar in spirit
to the one by Lipman in \cite{lipman})

Since  the local fundamental group is the fundamental group of
   $U -  \{P \}$ where $U$ is a good neighbourhood  of $P$  in $X$,
   Mumford considered the minimal normal crossings resolution
of the singularity, and derived the above theorem from the following.

Let $ D = \cup_i D_i $ be a  compact connected normal crossings divisor
on a smooth algebraic surface $S$, such that the intersection matrix
$(D_i  \cdot  D_j )$ is negative definite : then the local fundamental
group around $D$, i.e., the fundamental group $\Ga : = \pi_1 (T - D)$
   where $T$ is a good tubular neighbourhood of $D$, is trivial
if and only if $D$ is an exceptional divisor of the
first kind (i.e., $D$ is obtained by successive blowing ups
starting from  a smooth point
   of another algebraic surface).

Our purpose  here is threefold:

1) first, we want  to show that the theorem has more to do with
   a basic concept appearing in surface classification
rather  than with singularities; i.e., that the crucial hypothesis
is not that the matrix $(D_i  \cdot  D_j  )$ be negative definite, but
that the canonical divisor $K_S$ of $S$ be nef on $D$
   (this happens for a minimal model of a nonruled algebraic surface).
For the nonexpert: the condition that $K_S$  be nef on $D$
means that , if $g_i$  = genus of the smooth curve $D_i$,
then for each $i$ it holds: $2 g_i   - 2 \geq D_i^2 $.

As a matter of fact, this condition will only be needed for the curves
$D_i$  of genus zero, for which it reads out as  $  D_i^2 \leq -2$ .

\medskip
2) Second, since  the structure of the group $\pi_1 (D)$
is very well understood and there is an obvious surjection
   $\Ga = \pi_1 (T - D) \rightarrow \Pi : = \pi_1 (D)  $, we want to study in
general how big is the kernel $ \sK$
of this surjection.  Then  the  result is that under the above
nefness hypothesis each standard generator of $ \sK$, i.e.,
   each  simple loop $\gamma_j$ around a component $D_j$, is nontrivial in
   $\pi_1 (T - D)$.

More precisely, we would like to show that, outside of a well described
family of exceptions, this generator $\gamma_j$  has infinite order.

It is rather clear that, in order to have a very simple formulation,
   the hypothesis that $K_S$ be nef on $D$ is necessary.

In fact, if we let $D$ be a line in $\P^2$ , the local fundamental
group  around $D$ is trivial and we have $K_{\P^2} D= -3$;
   similarly happens if we take a $(-1)$ -curve (a smooth rational
   curve
with self intersection = -1, hence a curve with $K_S D = -1$).

A slightly more complicated example , obtained by blowing up
   the central point of a string of 4 $(-2)$-rational curves,
shows that the local fundamental group may be nontrivial,
yet some $\gamma_i$  may be trivial, if we do not use the
nefness assumption.

The simplest results we have in the direction explained above
are the following theorems A, B, C.

Among these , the following theorem A is,
as already said, the simplest one to be stated:

\begin{thm}
{\bf (Weak Plumbing Theorem A)}.

Let $ D = \cup_i D_i $ be a  connected compact (global)
normal crossings divisor
on a smooth complex surface $S$.

Assume further that the dual graph $ \sG$ of $D$ is a tree.

Let $\Sigma$ be the boundary of a good tubular neighbourhood
$T$ of $D$ , $ T = \cup_i T_i $.

The generator $\gamma_i$  of the kernel
$\cong \Z$ of
   $\pi_1 (T_i  - D_i ) \rightarrow \pi_1 (D_i )$
   has a nontrivial image in
   $\pi_1 (\Sigma ) \cong \pi_1 (T  - D ) $
if it holds true the stronger assumption  that
the canonical divisor $K_S$ of the surface $S$ is nef on
the components of $D$ of genus $0$, i.e., $ K_S D_i \geq 0$ for
each $i$ such that $D_i$   has  genus zero.
\end{thm}

\begin{rem}
    Please observe that we do not need $S$ to be compact:
this hypothesis would entail, by the Index Theorem, that
the positivity index of the matrix $( D_i  \cdot  D_j ) $ be $ \leq 1$.

Therefore, our result concerns all the 3-manifolds
$\Sigma$ which are boundaries of complex surfaces
obtained by plumbing smooth compact complex curves.
\end{rem}

More generally, holds the more precise

\begin{thm}
{\bf (Strong Plumbing Theorem B)}.

Let $ D = \cup_i D_i $ be a  connected compact (global)
normal crossings divisor
on a smooth complex surface $S$.

Assume further that the dual graph $\sG$ of $D$ is a tree.

Let $\Sigma$ be the boundary of a good tubular neighbourhood
$T$ of $D$ , $ T = \cup_i T_i $.
Then the generator $\gamma_i$   of the kernel $\cong \Z$ of
   $\pi_1 (T_i  - D_i ) \rightarrow \pi_1 (D_i )$ has
a nontrivial image in  $\Ga : = \pi_1 (\Sigma ) \cong \pi_1 (T  - D ) $
if

i)	 $D$ is minimal , i.e. , it is not obtained by blowing up a
(global) normal crossings divisor $D'$
and moreover either

ii-1)	after successively blowing down all the rational $(-1)$-
   curves we get a divisor $D'$ contained in a smooth
complex surface $S'$ and such that $K_{S'}$  is nef on the
components $D'_i$  corresponding to a $D_i$  of genus zero, or

ii-2)	 if $D_i$  has  genus zero, then its self intersection
is negative.
\end{thm}

3) Our motivation for studying these  questions came from
the study of
   topological characterizations of the existence
of fibrations on algebraic surfaces, especially in the noncompact
case, where (cf. \cite{cat00})  one has to
consider the fundamental group at infinity, which
is a disjoint union of  local fundamental
groups $\pi_1 (T - D)$.

   The goal is to  get  new and simpler
    variants of the characterizations of the Zariski open sets
which are the complement of a union of fibres of a fibration
    containing all the singular fibres.
These were given in  \cite{cat00}, theorem 5.7, for
constant moduli fibrations, and in \cite{cat03}, theorem 6.4, in the
general case.

Indeed, in these theorems there is one condition
   pertaining the fundamental group at infinity, namely
that, given a certain group homomorphism,
each $\gamma_i$  maps to a certain element of infinite order.

So, a natural question is: when does each $\gamma_i$
have infinite order in $\pi_1 (T-D)$?

We have some partial result concerning this question, which we hope 
to be able to
improve in the future

\begin{thm}
{\bf ( Plumbing Theorem C)}.

Let $ D = \cup_i D_i $ be a  connected compact (global)
normal crossings divisor
on a smooth complex surface $S$ satisfying the assumptions of the
previous Theorem A,
(we want again  for instance that the dual graph $\sG$ of $D$ is a tree).

Define $D$ to be {\bf elementary infinite} if either

1) $\sG$ is a linear tree and there is a curve of positive genus, or

2) $D$ is  a {\bf comb} (i.e., $\sG$ contains only one vertex of valency 3)
and there is a curve of positive genus, or all curves are of genus 0,
but we are not in the exceptional cases Va) and Vb).

Let $\Sigma$, $\Ga$, $\ga_i$ be as in the previous theorems:
then each $\gamma_i$
   has infinite order  in
   $\Ga  $
if there is a sequence of moves, consisting in successively removing 
curves $D_i$
which intersect two or more other curves, such that in the end
one is left with a bunch of disjoint elementary infinite pieces.

\end{thm}

Actually, since it can happen that the normal crossing
configuration be not minimal,  it would be certainly interesting to give
necessary and sufficient general conditions also for the nontriviality
of each $\ga _i$ (this might be very complicated, we fear).

For the  applications mentioned above, however, we need to treat the 
general case
and we may not restrict ourselves to the situation where the dual graph
is a tree, which is treated in this article.

As a matter of  fact, at some point  we thought we could easily
reduce the case where the dual graph is not a tree to the difficult case
where we have a tree: but about five years ago, when we were writing up
a first version of the article, we realized that this reduction argument
was not correct.

One reason why we want now to write down here  the tree case,
is because this article owes much to Guido Zappa.  When
I started to think  about these questions, I received a kind letter of Zappa,
which was somehow related to my election as a corresponding member
of the Accademia dei Lincei, and it was only natural to ask him
some question in combinatorial group theory.
Zappa not only answered,  providing a result which is included in the
article (cf. proposition \ref{zappa}), but he was very kind to continue
   to read and answer my letters.

Thus this article is particularly appropriate for this special volume of the
Rendiconti Lincei, dedicated to Guido Zappa. I am
 indebted to him,
   to  his wife Giuseppina Casadio and also to Antonio Rosati for 
orienting my choice
towards mathematics.  Giuseppina Casadio ran some afternoon seminars
in the Liceo Ginnasio 'Michelangelo' in the last year of  my
(classical studies) high-school.
There I learnt such basic things as, for instance,  congruences, and
I was encouraged to take part into the Mathesis competitions first
and  the mathematical Olympics later.
Rosati incited me over the summer to read parts of Courant
and Robbins'
book 'What is mathematics', and to apply for admission to the Scuola Normale
Superiore di Pisa.

In Pisa the education was very analysis oriented, but later on in my
life I discovered
in myself something of an algebraist's soul which was longing to learn more.

For this part of my soul Zappa was a reference
figure, and
I was later quite happy to have finally a chance,
during the Meetings of the Accademia, to discuss mathematical questions
with him.

Another reason to write this article now is to take up the problem 
again, with the
hope of finding soon the solution to the general case, and, even 
more, to propose
the further investigation
of these three-manifolds fundamental groups.

For instance,  other general
   interesting questions are  in our opinion:

1) how big is the kernel $ \sK$ of $\pi_1 (T\setminus D) \ra \pi_1 (D)$ ?

2) What properties does $ \sK$ enjoy, when is it for instance not
finitely generated (cf. \cite{cat03}, definition 3.1 and lemma 3.4)?

\section{A presentation of the local fundamental group}
\label{pres}

Let us  first of all set up the notation for our problem.

We have $S$ a smooth complex surface, and a compact connected
global normal crossings divisor  $ D = \cup_i D_i $
contained in $S$,
thus each $D_i$  is a smooth curve of genus $g_i$ and has a
good tubular neighbourhood $T_i$ which is a 2-disk bundle over
$D_i$.

$T_i  \setminus D_i $  is homotopically equivalent to its boundary
$\Sigma_i$  , which is an $S^1$ -bundle  over the compact
Riemann surface $D_i$, and is completely classified by
   its Chern class, i.e., by the self-intersection number of $D_i$ in
$S$, as we are going to briefly recall.

Let us denote by $m_i$  the opposite  of  the self intersection
number of $D_i$ , so that we have $D_i^2= - m_i$.

Let now $q$ be a point of $D_i$ : then the bundle
$\Sigma_i \rightarrow  D_i$  is trivial over $D_i \setminus \{q\}$, and also
   over a neighbourhood $V$ of $q$ . The respective trivializations
are clear if we identify topologically  the associated line bundle as the 
line bundle corresponding to the divisor  $ - m_i \ q$.

Since $(D_i -q) \cap V$ is homotopically equivalent to $S^1$,
and the glueing map on $S^1 \times S^1$ reads out
( we choose  the trivialization over $D_i \setminus \{q\}$ in the source, and the 
one over  $V$  in the target)

$$(z,w) \rightarrow (z, z^{- m_i} w),$$  from the I van Kampen
Theorem (cf. e.g. \cite{dR}) we derive a presentation for the 
fundamental group of
$\Sigma_i$, which determines the central extension

$$1  \rightarrow \Z  \gamma_i \rightarrow \pi_1 (\Sigma_i)
   \rightarrow \pi_1  (D_i ) \rightarrow 1$$

provided by the homotopy exact sequence of the $S^1$ -bundle.

   In fact, in the inverse image of $D_i \setminus \{q\}$, $\cong 
D_i \setminus \{q\}
\times  S^1$ we take the lifts of some standard generators of the free group
$\pi_1 (D_i -q)$, we use for these lifts the usual notation
$a_1 (i) , b_1 (i) ,  ...  a_{g_i} (i) , b_{g_i}(i)$
   (recall that $g_i$  is the genus of $D_i$), and moreover
we let $\gamma_i $   be the generator of the fundamental group
of the fibre $S^1$, with the standard complex
counterclockwise orientation.

Since the fundamental group of a Cartesian product
   is a direct product, it follows, as already mentioned,
   that $\gamma_i $  commutes with
all other generators.

  From the glueing map we get the single further relation :

$$ \prod _{h=1,...g_i}  [a_h (i) , b_h (i) ]=  \gamma_i^{- m_i}.$$

If we take now a good tubular neighbourhood $T$ of $D$ which is
the union of the $T_i$` s, we may assume moreover
( by shrinking the $T_i$` s,  and by the implicit function theorem),
   that the intersection $T_i \cap T_j$ be biholomorphic to

$$\{ (z_1 ,z_2 ) |  \ |z_1 z_2| \leq 1, |z_i| \leq 2\},$$

   where $z_1  = 0, \ z_2 = 0,$ are the respective local equations of
$T_i  , T_j$,   at the point $p_{ij} :=D_i  \cap D_j$.

In each $D_i$  let us consider a path $L_i$ homeomorphic to
a segment  and going through
   all the points $p_{ij}$  and let us mark a point $q_i  \in L_i$
    different from  all the $p_{ij}$ ' s.

We may easily assume that we get thus a linear tree $L_i$
with the above points as vertices.

Set $ L =\cup_i L_i$  , thus $L$ is naturally a graph.

It is important to notice that $\Sigma$ has a natural projection
onto $D$, such that outside the points $p_{ij}$  we have a
fibre bundle with fibre $S^1$, whereas the fibre over $p_{ij}$  is
   $ \cong  S^1 \times S^1$.

In fact,  the local picture is given by

$$ T_i \cap T_j  = \{ (z_1 ,z_2 ) |  |z_1 z_2| \leq 1, |z_i|
\leq 2\},$$ thus locally
$$\Sigma = \{ (z_1 ,z_2 ) |  |z_1 z_2| = 1, |z_i|
\leq 2\}  \cong S^1\times  S^1 \times [1/2, 2] ,$$

where the homeomorphism is given by the map sending $(z_1 ,z_2)$
to

   $(z_1/|z_1| , z_2/|z_2| , |z_1| ).$

The projection sends   $S^1\times  S^1 \times \{1\}$  to
$(0,0)$, whereas e.g.  the observation that

$S^1\times  S^1 \times [1/2, 1)$ is an  $S^1$-bundle over
$S^1 \times [1/2, 1) \cong$ punctured disk in the  $z_2$ plane,
allows to define the projection for  $|z_2| \geq 1 $ as sending
$(z_1 , z_2) \rightarrow ( 0, z_2  (|z_2| - 1))$ ,
and symmetrically for
$|z_1| \geq 1 $.

It is quite easy to see then that we can find a section of
$ \Sigma | _L \rightarrow |L$, so we think of L as $
\subset \Sigma | _L $.

Since the restriction of the fibration $ \Sigma_i  \rightarrow
D_i$  to  $L_i$  is trivial , we obtain that,
   up to homotopical equivalence, $ \Sigma | _L \rightarrow L$
   is obtained from the manifolds $L^0_i  \times S^1 $
($L^0_i$  being a tubular neighbourhood of $L_i$  in $D_i$)
as follows.

We replace the product $B^2_{ij} \times S^1$ ($B^2_{ij}$
   being an open  2-dimensional
ball around $p_{ij}$  in $D_i$) by a product $A^2_{ij} \times S^1$
(($A^2_{ij}$
   being a  2-dimensional annulus around $p_{ij}$  in $D_i$,
   $A^2_{ij} \cong  S^1 \times [1/2, 1)$ ).

Then we glue together the pieces $A^2_{ij} \times S^1$ and
$A^2_{ji} \times S^1$  identifying the (inner)
   boundaries $S^1 \times  S^1$ .

We make now another arbitrary choice for our presentation,
namely, since the graph $L$ is connected , we may take a
connected subtree $L' \subset L$ containing all the points $q_i$.

We let one of them, say $q_0$ , be the base point : for each $q_i$
    we get a canonical path in $ L'$ from $q_0$ to $q_i$ ,
whence a canonical basis of $ \pi_1 (L)$ is given by the loops
$\lambda_{ij}$  , for $p_{ij}$ not $ \in  L'$ , obtained going from
$q_0$ to $q_i$ along the canonical path, then going to $p_{ij}$
inside $L_i$  ,
   then to $q_j$   inside $L_j$,  then back to $q_0$
again along the canonical path.

The above description makes it clear that , exchanging the role
of the two indices $i,j$, we get $\lambda_{ji} = \lambda_{ij}^{-1}$.

Let $\gamma_i$  be the positively oriented generator of the
infinite
cyclic fundamental group of $(L^0_i  \times S^1 ) \cup L'$:
then we find immediately the following presentation for the
fundamental group of $\Sigma$ restricted to $L^0$ ( $L^0= \cup_i
L^0_i$).

\begin{itemize}
\item
(2.12) {\bf Generators}:
\item $\gamma_i$  , for each $i$,
\item
$\lambda_{ij} $  , for $p_{ij}$ not $ \in  L'$.

In order to get the relations , set , for each $p_{ij} \in L$,

\item
				$\gamma_{ij} = \gamma_j$, for $p_{ij}
\in L'$, and
\item
			     	$\gamma_{ij} = \lambda_{ij} \gamma_j
\lambda_{ij}^{-1}$,
   for  $p_{ij}$  not in  $L' $,

with the above convention that $\lambda_{ji} = \lambda_{ij}^{-1}$.

Then we get the
\item
(2.13){\bf Local Commutation Relations}:
				$[\gamma_i , \gamma_{ij}] = 1 $ (for
each $p_{ij} \in L$).
\end{itemize}

To complete the presentation of $ \pi_1 (\Sigma)$,
we use several times again the First van Kampen theorem
( cf. [dR]), adding $ \Sigma |_{(D_i  - L_i  )}$ to  $\Sigma$
   restricted to $L^0$ .
   Note that the $S^1$-bundle $ \Sigma_i  \rightarrow D_i$
is trivial on $L^0_i$  , and also on $D_i  - L_i$.

The corresponding fundamental group is obtained as
amalgamation by $ \Z \gamma_i  \times \Z \mu_i$  of the
free product of  the following two groups: the direct product
$\F_{2g_i} \times  \Z \gamma_i$ ($\F_{2g_i} $ = free group in $2 g_i $
generators) and the cyclic group $\Z \gamma_i$.

Here, $\mu_i$   maps on the one side to the standard relation
for the fundamental group $ \Pi_{g_i}$   of a compact curve
of genus $g_i$, on the other side it maps  to $\gamma_i^{m_i}$.

Now, $\mu_i$   is no longer trivial in $ \pi_1 (L^0 )$,
so we get the following extra

\begin{itemize}
\item
(2.14) {\bf Generators} : $a_1 (i) , b_1 (i) ,  ...    a_{g_i} (i),
b_{g_i} (i)$ , for each $i$,
\item
(2.15) {\bf Main relations}:
$$ \prod _{h=1,...g_i} [a_h (i) , b_h (i) ]=  \gamma_i^{- m_i}
\prod _{j} \gamma_{ij} .$$

Moreover, since we have a direct product $\F_{2g_i} \times  \Z \gamma_i$,
   we should not forget the obvious relations :
\item
(2.16) {\bf Global Commutation relations} :
$[a_h (i) , \gamma_i] = [\gamma_i , b_h (i) ] = 1 .$

\end{itemize}

\section{Presentation of a simplified group }
\label{group}

Summarizing the result of the previous section,
we have gotten the following finitely presented group $\Gamma$  with :

\medskip

{\bf GENERATORS}:

   $\gamma_i$  , for each $i$,

$\lambda_{ij} $  , for $p_{ij}$ not $ \in  L'$,

$a_1 (i) , b_1 (i) ,  ...    a_{g_i} (i),
b_{g_i} (i)$, for each $i$.

\medskip

{ \bf RELATIONS} :

\begin{itemize}
\item
$[a_h (i) , \gamma_i] = [\gamma_i , b_h (i) ] = 1 ,$ for each $i,h$
\item
$ \prod _{h=1,...g_i} [a_h (i) , b_h (i) ]=  \gamma_i^{- m_i}
\prod _{j} \gamma_{ij} $ for each $i$,
\item
$[\gamma_i , \gamma_{ij}] = 1 $ (for each $p_{ij} \in L$) and where

	I)			$\gamma_{ij} = \gamma_j$, for $p_{ij} \in L'$,

	II)		     	$\gamma_{ij} = \lambda_{ij} \gamma_j
\lambda_{ij}^{-1}$,
   for  $p_{ij}$  not in  $L' $, and recall also

III) $\lambda_{ji} =
\lambda_{ij}^{-1}.$

\end{itemize}

\begin{rem}
The projection $ p : \Sigma \ra D$ induces a surjection of fundamental groups
$ \Ga \ra \pi_1(D)$ with kernel $\sK$ normally generated by the $\ga_i$'s.
In fact, setting in the above presentation $ \ga_i =1 \ \forall i$, we get
a free product of the fundamental groups $\pi_1(D_i)$ with  the free group
generated by the $\la_{ij}$'s (observe that  $\lambda_{ji} =
\lambda_{ij}^{-1}$, whence the rank of this free group is equal to the
first Betti number of $L$).
\end{rem}

\begin{defn}

The associated  {\bf simplified } finitely presented {\bf group }
$\Gamma'$
    is the following group $\Gamma'$    with :

GENERATORS :

   $\gamma_i$  , for each $i$,

$\lambda_{ij} $  , for $p_{ij}$ not $ \in  L'$,

$a_i , b_i$, for each $i$ such that $g_i \geq 1$.

RELATIONS :

( Global commutation relations)
   $[a_i , \gamma_i] = [\gamma_i , b_i] = 1 ,$ for each $i$

( Main relations)
$ [a_i , b_i]=  \gamma_i^{- m_i}
\prod _{j} \gamma_{ij} $ for each $i$,

(Local commutation relations)
$[\gamma_i , \gamma_{ij}] = 1 $ (for each $p_{ij} \in L$) where,
as above,
   $\gamma_{ij} = \gamma_j$, for $p_{ij} \in L'$, else
( keeping in mind :
$\lambda_{ji} =
\lambda_{ij}^{-1}$) $\gamma_{ij} =
\lambda_{ij} \gamma_j \lambda_{ij}^{-1}$.
\end{defn}

\begin{rem}\label{simplified}
We can restrict ourselves to prove our results for the
simplified groups $\Gamma'$ , which are also obtained from a
plumbing procedure, replacing the (smooth) curves of genus
$ \geq 2$
by genus $1$ curves.

In fact, the simplified  group $\Gamma'$  is a homomorphic
image of $\Gamma$  , being obtained by imposing the further
   relations

$a_h (i) =  b_h (i) = 1$ , for $ h \geq 2 $.

Thus, if $\gamma_j$   is nontrivial,
respectively of infinite order, in the simplified group $\Gamma'$
it is so a fortiori in the group $\Gamma$ .
Moreover, observe that our hypotheses only concern the nullity
or positivity of  the genus of $D_j$  , and not its precise value.

For instance, the minimality of D in the category of normal
crossing divisors amounts to the nonexistence of rational curves
with self intersection $= -1$, and meeting at most two other
curves each in at most one point.
Thus, we see easily that  the hypothesis i) of B)
is still verified for the simplified group,
likewise for the hypothesis of A).

   We may have however that the canonical divisor $K'$ of the
simplified  surface could not be nef, since if  there is a
component $D_i$ with genus $ \geq 2$ , in the new configuration
$C$ we get a corresponding $C_i$ with genus $1$ and
$ K'  C_i = -  C_i^2 = - D_i^2 = -  (2 g(D_i )-2) + K  D_i \ $,
which may become negative.
\end{rem}

The proof of the main theorems follows by a reduction
step which we  examine in the next section.

\bigskip

\section{Reduction to the case of a graph of rational curves}
\label{ rat graph red}

Recall  that we are working in the simplified group.

In the case where we get a component of genus 1, we will be able to
simultaneously remove the generators $a_j, b_j$, and replace the number $m_j$
by  any arbitrary integer $n_j$  (in fact, one could say that we can have
$n_j = \infty $, meaning that the corresponding main relation disappears).

If we can achieve this, certainly the nefness condition on the new
configuration will continue to hold. To this purpose, let us  fix the
index  $j$,  let us write
$$   a_j : = a, b_j  : =b ,\ga : =\ga_j ,$$    and let us consider
the group G generated
by generators
\begin{itemize}
\item
$\ga_i$  , for each $i$,
\item
$a_i , b_i $, for the $i$ 's such that  $g_i  \geq 1$, and $i \neq j$

and by relations :
\item
$[a_i  , \ga_i  ] = [\ga_i  , b_i ] = 1 $,  for each $i \neq j$
\item
$[a_i  , b_i  ] =   \ga_i ^{-m_i} \prod_h \ga_{i h} $, for each $i \neq j$
\item
$[\gamma_i , \gamma_{ih}] = 1 $ (for each $p_{ih} \in L$).
\end{itemize}
The group $\Ga$ is obtained from $G$ by adding generators $a,b$,  and relations
\begin{itemize}
\item
$[a , \ga ] = [\ga , b] = 1 $ where  $\ga : =\ga_j $  is an element of G ,
\item
$ [a, b]=  \ga^{-m} \prod_h \ga_{jh} $.

We may rewrite the last relation simply as
\item
$ [a, b]=  \ga '' $.
\end{itemize}
Note that , in the group $G$ , $[\ga , \ga '' ] = 1$ , since $\ga$
commutes with each
$\ga_{jh} $.

We use now:

\begin{prop}\label{zappa}
    Given a group $G$ , and elements , $\ga , \ga '' \in G$
such that  $[\ga , \ga '' ] = 1$, let $\Ga$ be the group obtained as
the quotient of
the free product of $G$ with a free group generated by two generators $a,b$,
by imposing the following relations :
$$ [a , \ga ] = [\ga , b] = 1 , [a, b]
= \ga ''. $$Then the natural homomorphism of  $G$ into $\Ga$ is injective.
\end{prop}

\Proof. We
consider the quotient group $\Delta$  of $\Ga$ obtained by adding the
commutation
relations $[a , \ga '' ] = [\ga '' , b] = 1$ . An equivalent way to describe
$\Delta$ is
the following.

Let $H$ be the Heisenberg group generated by generators
$a,b,c$ and with relations  $[a,b] = c , [a,c] = [b,c] = 1$. $H$ is a two step
nilpotent group with infinite cyclic centre generated by $c$ , and
abelianization free of rank 2 . The elements in $H$ can be uniquely written
as words  $a^m b^n c^k$, where $k,m,n$ are integers.

Then we can define $\Delta$  as
the quotient of the free product of $H$ and $G$ , modulo the relations
$$\ga ''= c
, [a , \ga ] = [\ga , b] = 1. $$
At this point we are not able to have a unique
representation for the elements of $\Delta$ , but we follow an idea
of Guido Zappa.

Namely, we observe that every element of  $\Delta$ can be written as a product

$$h = g_0  a^{m(1)} b^{n(1)} g_1  a^{m(2)} b^{n(2)}
    \dots g_{r-1} a^{m(r)} b^{n(r)}  g_r  , $$
where
each pair of exponents $(m(j),n(j))$ is $ \neq (0,0)$, $ g_0  , \dots
g_r $  are elements
of $G$ and we can assume that $ g_1  , \dots  g_{r-1}$  do not belong to the
subgroup $B$ generated by  $ \ga , \ga " $ in $G$.
   ( whereas ,  $g_0$  and $g_r $ could be even trivial).

There remains to see
when two such products yield the same element $h$ . Notice that the
condition (*) that $ g_1  , \dots  g_{r-1}$  do not belong to $B$
follows from the
property that $r$ be minimal.

We claim that $r$ is uniquely determined , and
that the  only allowed transformations of the minimal  representation
are obtained by letting
   factors $\ga , \ga " $ commute with $a$, resp. $b$.

More precisely, we claim
that we get an equivalent minimal product iff :

\begin{itemize}
\item
   we replace each
respective element $g_i$  ($ =g_1  , .. $ or $ g_{r-1}$) multiplying
it by an element
$ g \in B $, and correspondingly :
\item
   if $g_i$  is replaced by $g_i  g$ , then
$g_{i+1}$  is replaced by $ g^{-1} g_{i+1}$ ,
\item
   if $g_i$  is replaced by $g g_i  $  , then
$g_{i-1}$  is replaced by $ g_{i-1} g^{-1} $

\end{itemize}

This means that, for each $i$ , the
exponents  $(m(j),n(j))$ are uniquely determined ; moreover, the double
coset  $B g_i  B$ is uniquely determined, and finally the product
$g_0   \cdots g_r $
is uniquely determined . In particular, it follows that our element is in $G$
iff  $r=0$, and in this case the representation is unique , what is
precisely the
assertion of the proposition .

To establish our claim , let us consider the
equivalence classes of the products $h$ described above. It suffices to show
that we have an action of the generators of the group $\Delta$, which satisfies
the defining relations for $\Delta$.  This is clear for the elements
of the group $G$,
and also for the generators $a,b$, and an easy verification show that the
relations are satisfied.

				\qed

\begin{rem}\label{arbitrary}

Notice that , if we fix an integer $n_j$  and in the group $G$ we add the
relation
$$ 1=   \ga_j ^{-n_j} \prod_h \ga_{i h },$$

we have the corresponding fundamental group of the graph of curves where
the elliptic curve $C_j$  with self intersection ($-m_j$) has been replaced by
a smooth curve $\cong \PP^1$ with self intersection ($- n_j$). We can therefore
by induction reduce to the  case of  a graph of rational curves .
\end{rem}

\section{The case of a tree of smooth rational curves}
\label{ rat tree }

We have here a presentation with

{\bf GENERATORS}:

   $\gamma_i$  , for each $i$,

\medskip

{ \bf RELATIONS} :

\begin{itemize}

\item
$ 1 =  \gamma_i^{- m_i}
\prod _{j} \gamma_{ij} $ for each $i$,
\item
$[\gamma_i , \gamma_{j}] = 1 $ (for each $p_{ij} \in L$)

\end{itemize}

We would like first
to show the necessity of the nefness hypothesis in Theorem A.

\begin{ex}\label{blowup}  Consider a diagram of type $A_n$, i.e., a linear tree
with
$n$ vertices.

Then our group, as we shall shortly see,  is generated by :
$\ga_1  , \dots\ga_n$ , with relations
$$\ga_1 ^2 =  \ga_2,   \ga_2 ^2 =  \ga_ 1\ga_3,   \ga_3 ^2 =
\ga_2 \ga_ 4,   \dots
   \ga_{n-1} ^2 =  \ga_{n-2 }\ga_n,   \ga_n ^2 =  \ga_{n-1}.  $$
    Therefore, the group is cyclic ,
generated by   $\ga: = \ga_{1}$,  with
   $ \ga_{1} ^{n+1} = 1,$ and we have  $\ga_{i}  =  \ga_{1} ^i  .$

   Let  $n= 4$, and let us now blow up the
central point of intersection between   $C_2$  and $C_3$.

We obtain then a
new generator $\ga '$ (the loop around the exceptional curve) and the relation
$\ga ' =  \ga_2 \cdot \ga_3$, but then $\ga ' = \ga_2 \cdot \ga_3 = 1$ !

\end{ex}

We have to recall, in the case where we have a tree of rational curves on a
complex surface , that the condition that the divisor $K_S$ is nef reads out
as

1) $D_i  ^2 \leq  - 2 . $

      If we are on an algebraic surface, the index theorem says that

2) the intersection matrix $(D_i  \cdot  D_j  )$ has positivity index
$b^+ \leq 1 . $

An easy example where 1) holds but $b^+ = 1  $ is provided by a  tree of
rational  $(-2)$ curves , where all curves meet a central one (the dual
graph is a star ).

In fact, then , if  $D_0$ is the central curve , we have
$$(m D_0  + D_1  + \dots  D_n  )^2 = 2 ( -m^2 +mn - n ),$$ which is
positive for
$ 1 <
m < n -1 .$

   Then the group is generated by $\ga_1  , \dots \ga_n , \de $, with
relations
$$ \ga_i ^2 =  \de,   \de ^2 =  \ga_1 \cdot \ga_2\ \cdots \ga_n. $$
In this case the
Abelianization is the direct sum of cyclic groups of respective orders
$2 (n-4) , 2, \dots  2$ , with generators induced by the respective
residue classes of
$\ga_1  ,\ga_1 ^{-1} \ga_2,  \dots \ga_1 ^{-1} \ga_{n-1} $,
whence
here our standard generators have even a nontrivial image in the maximal
Abelian quotient.

We proceed now to analyse the different cases.

{\bf5 A :}   CASE OF A LINEAR TREE OF RATIONAL CURVES

\begin{lem}\label{ltree}
   Assume that we have a linear tree of $n$  smooth rational curves with self
intersection ($-m_i$ ) , where  $m_i \geq 2$ .

Then , setting inductively
$ a_1 : = 1, a_2 : = m_1 ,  a_{i+1} : = m_i \cdot  a_i - a_{i-1},$ then

   1) $ a_{i+1}  > a_i $;

2)  our
group $\Ga$ is a cyclic group of order  $a_{n+1}$ , generated by $\ga_1$;

3) the element  $\ga_i$  equals $\ga_1 ^{a_i} $, and is not trivial .
\end{lem}
   \Proof
We can write our relations among  $\ga_1,\ga_2,  \dots\ga_n$  as
$$ \ga_1^{m_1} = \ga_2, \ga_2^{m_2} = \ga_1 \ga_3,   \dots
\ga_i^{m_i} = \ga_{i-1} \ga_{i+1}, \dots  \ga_{n-1}^{m_{n-1}} =
\ga_{n-2} \ga_n,
\ga_n^{m_n} = \ga_{n-1} . $$
We easily obtain then
     $$\ga_{i+1}=\ga_{i-1}^{-1}  \ga_i^{m_i} = \ga_1 ^{-a_{i-1}} \cdot
   \ga_i^{a_i m_i} = \ga_1 ^{a_{i +1}}  ,$$ which proves the first
part of assertion 3) , and the last relation on the other hand yields
$\ga_1^{a_{n +1}} = 1,$
which proves assertion 2 .

   Notice that
$$a_{i+1} - a_i =  m_i \cdot a_i - a_{i-1}- a_i  =
(m_i -1 )\cdot  a_i - a_{i-1} > 0$$  since $ m_i \geq 2 $ and since
by induction $ a_i > a_{i-1}$.

Whence, assertion 1) is proved, and simultaneously we have shown
that each $\ga_i$  is not
trivial .

									\qed

\begin{rem}\label{blowdown}
The proof of the above lemma shows that in
any case the local fundamental group of a tree of rational curves is
cyclic, of order $a_{n +1}$ , if $a_{n +1}$ is nonzero.

Assume now that all the
numbers $m_i$ are strictly positive. Then, if $m_i = 1$, we obtain
$ \ga_i = \ga_{i-1}
\ga_{i+1},$ and since the group is abelian , we may rewrite the relation
$\ga_{i-1} ^{m_{(i-1)}} = \ga_{i-2} \ga_i$ as  $\ga_{i-1} ^{m_{(i-1)} -1} =
\ga_{i-2}  \ga_{i+1}  $
and similarly
$\ga_{i+1} ^{m_{(i+1)}} = \ga_i \ga_{i+2} $ becomes
$\ga_{i+1} ^{m_{(i+1)} -1} = \ga_{i-1} \ga_{i+2}$ .

   This has the
obvious geometrical meaning that we can blow down all the (-1) curves,
and then if at the end of the process $K$ remains nef, our remaining
elements $\ga_i$ are not trivial.
\end{rem}

\begin{rem}\label{infty}
Assume that we let $ m_i \ra \infty$. Then also  $ a_{i+1} \ra \infty$,
hence $ a_{n+1} \ra \infty$, whereas $ a_j$ remains constant for
$ j \leq i$. Hence, $ ord (\ga_j)  \ra \infty$ for
$ j \leq i$. Changing the linear order
of the linear tree to its inverse, we see that $ ord (\ga_j)  \ra \infty$ 
also for
$ j \geq i$. 
\end{rem}

{\bf 5 B : }  REDUCTION TO THE CASE OF A COMB OF RATIONAL CURVES

\begin{lem}\label{G1G2}
   Let  $G_1 $, $G_2 $ be groups and let  $a_i$ be nontrivial 
elements in $G_i$,
for  $i= 1,2 $, such that moreover $a_2$ has infinite order in $G_2 $ .

If $\Ga$ is the
quotient of the free product  $G_1  * G_2 $  by the relation  $a_1
\cdot a_2 = 1$ , then
the natural homomorphism of  $G_1 $ in $\Ga$ is injective. Moreover,
if $a_1$  does
not generate $G_1 $ and $a_2$  does not generate $G_2 $, then $\Ga$
is always an infinite
group.

\end{lem}
\Proof  . The desired claim follows if we show that the elements in
$\Ga$  are represented by elements of the set $\sW$  of equivalence classes of
'good' words
$$w = g_1(1) \cdot g_2(1) \cdot g_1(2) \cdots g_1(k) \cdot g_2(k)
\cdot  g_1(k+1) ,$$
where  $g_2(i)$
does not belong to the subgroup generated by  $a_2$ , for $1 \leq i \leq k $, 
and $g_1(j)$
does not belong to the subgroup generated by  $a_1$ , for $2 \leq j \leq k $,

and  $w$  is equivalent to $w'$ if and only if  the following conditions hold:

1) k = k'

2) there exist
integers ( "$r$"for right, "$\la$" for left )  $r_1 , \la_2 , r_2 , \la_3 ,
\dots  r_k , \la_{k+1}$,  such that the word $w'$ equals

$$ (g_1(1) a_1^{r_1} )\cdot (a_2^{r_1}  g_2(1) a_2^{\la_2} )
(a_1^{\la_2} g_1(2) a_1^{r_2} )
\cdots  (a_1^{\la_k}
g_1(k) a_1^{r_k}  ) \cdot (a_2^{r_k}  g_2(k) a_2^{\la_{k+1 }}) \cdot
(a_1^{\la_{k+1 }}
    g_1(k+1) ).$$

We let the elements of $\Ga$ operate by left multiplication as follows :
\begin{itemize}
\item
for
$\ga_1 \in G_1 $ we let  $\ga_1 w  : = (\ga_1 g_1(1) )\cdot g_2(1) \cdot g_1(2)
   \cdots g_1(k) \cdot
g_2(k) \cdot g_1(k+1)$,
\item
for $\ga_2 \in G_2 $  not in the subgroup generated by $a_2$ we
let  $$\ga_2 w : =  e_1 \cdot \ga_2 \cdot g_1(1) \cdot g_2(1) \cdot g_1(2)
   \cdots g_1(k) \cdot
g_2(k) \cdot g_1(k+1),$$ ( $e_i$ being the identity element of $G_i$)
, while  we set
\item
   $ a_2^r  w
: = a_1^{-r } w .$

\end{itemize}

   We obtain a homomorphism of each $G_i$ into the group $\sS (\sW)$  of
permutations of $(\sW)$, and moreover the transformation associated
to $a_1\cdot  a_2$
is by definition the identity, whence we get a homomorphism of $\Ga$ into
$\sS (\sW)$.

Moreover , $\Ga$ acts transitively on $\sW$ . Representing each element of
$\Ga$ by a good word $w $, we see that  if $w$ is the identity
this  implies that $k= 0$ ,
and $g_1(1) = e_1$.

Thus the action on $e_1$ establishes a bijection between $\Ga$
and $\sW$,  in particular since the words with $k=0$ correspond to the elements
of $G_1$ , $G_1$ injects into $\sW$, whence into $\Ga$. Notice
finally that if $a_2 $
generates $G_2$ then $G_1$ is isomorphic to $\Ga$, similarly if
$a_1$  generates $G_1$.

Whereas, if $a_i $ does not generate $G_i$ ,then $k$ can be arbitrarily high,
whence $\Ga$ is surely infinite.

		\qed

\begin{cor}\label{3ormore}
   Let  $G_1 ,  \dots  G_r$ be groups and let  $a_i $, for $i= 1,... r
$, be a nontrivial
element in $G_i $. If  $\Ga$ is the quotient of the free product
$ G_1 * G_2 * \dots
* G_r $ by the relation $a_1\cdot  a_2 \cdots a_r = 1$ , then, for $r
\geq 3$, the natural
homomorphism of $G_1$ in $\Ga$ is injective. Moreover, if  $r \geq 4$, then the
group $\Ga$ is infinite.
\end{cor}
   \Proof.  Apply lemma \ref{G1G2} , considering that $ a_2 \cdots
   a_r$ is an element of infinite order in $  G_2 * \dots
* G_r $  .  In the case $r \geq 4$, apply the lemma to
   $G_1 * G_2$  and $  G_3 * \dots
* G_r $ , taking into
consideration that both are infinite and not cyclic.

		\qed

With the aid of the foregoing corollary we are able to reduce the proof
of our main results to a very special case.

\begin{prop}\label{combred}
   Let $\ga_i$  be one of our generators of the group $\Ga$, in the 
case where the
hypotheses of theorem B are satisfied: then $\ga_i$ is nontrivial 
except possibly
if the tree is nonlinear and the curve $D_i$  is the only one which 
intersects at
least three other irreducible components of $D$ (we shall then say 
that the tree
is a {\bf comb}, and that $D_i$  is the rim of the comb).
\end{prop}
\Proof.
   The case where
the tree is linear was already dealt with .

So, let us assume that there
exists a curve $D_j$ , with $ i \neq j$ such that  $D_j$ intersects
at least three
other irreducible components of $D$ . Let us consider the group $G$ obtained
as the quotient of $\Ga$  gotten by setting $\ga_j  = 1$.

If  $D - D_j$ (this denotes the difference as divisors, and not as 
sets)  has $r$
connected components $D(1) ,    \dots  D(r)$ , we see immediately that $G$ is
the quotient of the free product $G_1 * G_2 *  \dots  * G_r $ by the
relation $a_1\cdot
a_2 \cdots a_r = 1 $, where $G_h$ is the fundamental group of  the boundary of
a good tubular neighbourhood of $D(h)$ , and $a_h$ is the loop around the
unique irreducible component of $D(h)$ meeting $D_j$  . By our corollary , and
since by induction we may assume that each $a_i$ , $i= 1, \dots r $, 
is nontrivial, we obtain that each $G_h $ injects into $G$ , and a fortiori
into $\Ga$.

Whence , all elements $\ga_i$ with $ i \neq j$ are nontrivial.

				\qed

{\bf  5 C : }  THE RIM  OF A COMB OF RATIONAL CURVES.

\medskip

Assume that we have a unique curve $D_j$  such that $D - D_j$
    has $r \geq 3$
connected components $D(1) ,    \dots  D(r)$ ,  each being a chain of
smooth rational
curves. Set for convenience $\ga : = \ga_j $ .

We shall then say as before that we
have a COMB with RIM $D_j$ and with  STRINGS $D(1) ,    \dots  D(r)$.

Then ,
for each chain $D(h) $, we can order the generators in such a way that we
obtain relations
$$\ga_1 ^{m_1} = \ga_2, \ga_2 ^{m_2} = \ga_1 \cdot \ga_3, \dots
\ga_i ^{m_i} = \ga_{i-1} \cdot \ga_{i+1},
\ga_{n-1}^{m_{(n-1)} }= \ga_{n-2} \cdot  \ga_n .$$

   Proceeding as in section 5A ) , we infer that $\ga
=\ga_1^{a_n}$ , where $a_n > 0$ is defined inductively as in 5A).

Finally , letting
$(-m) $ be the self intersection of $D_j $ , we obtain a relation
$$\ga^m = \beta_1 ^{d_1} \cdot
\beta_2 ^{d_2}  \cdots \beta_r ^{d_r}  ,$$
   where the $\beta_h$ 's are the loops, for each chain $D(h) $,
around the  end opposite to  $D_j  $.

We are left with the following

\begin{thm}\label{rim}  Let $\Ga ( m, b_1 , b_2 , \dots b_r ; d_1 
,d_2 , \dots d_r
),
$  for
integers $ m \geq 2$ , $b_i > d_i \geq  1$, be the group generated by

i) generators  $\ga, \beta_1 ,\beta_2 , \dots \beta_r $, and relations

ii) $\ga = \beta_1 ^{b_1}=
\beta_2 ^{b_2} =  \cdots \beta_r ^{b_r}  $ (
recall that the integers $b_h$ are $\geq 2$ ), and

iii) $\ga^m = \beta_1 ^{d_1} \cdot
\beta_2 ^{d_2}  \cdots \beta_r ^{d_r}  .$

Then the ( central ) element $\ga$ is nontrivial inside $\Ga$ and
indeed of infinite order unless we are in the following exceptional
cases with  $r=3$, and where $ c= 1,2$, and $ 1 \leq t \leq n-1$:

Va ) $(b_1 , b_2 , b_3 ) = (2,2,n)$ , $ n \geq 2$,
$(d_1 , d_2 , d_3 ) =(1,1,t)$

   Vb ) $(b_1 , b_2 , b_3 ) =  ( 2,3, n) $, $3  \leq n \leq 5$,$(d_1 ,
d_2 , d_3 ) =
(1,c,t)$.

\end{thm}

\Proof

{\bf  Step I .}

   We may assume that G.C.D. $(b_i , d_i ) =1$ for each $i$.

This
is a consequence of the following Logical Principle Lemma of
Combinatorial Group Thery.

\begin{lem}\label{logical}{\bf (Logical Principle Lemma)}

   Let $G$ be a finitely presented group
$$G = \langle  \beta_1 ,\beta_2 , \dots \beta_r  |  R_1 ( \beta)=  \dots
R_s ( \beta )=  1 \rangle. $$
Then , setting  $\beta_1= \beta^ k $, i.e. , taking the new
group $ G'' : = G *\Z /   \langle \langle\beta_1  \beta^{ -k} \rangle \rangle$, 
we get  $ ord
_{G'' } (\beta) =
k \cdot ord_G  (\beta_1)$,
while, for $ j \geq  2$, $ ord _{G'' } (\beta_j) = ord_G  (\beta_j)$.
\end{lem}

\Proof
The situation is a particular case of lemma \ref{G1G2},
with  $a_1 =\beta_1 $,  and with $a_2 = \beta^{ -k }$.

The injectivity of the map $G \ra  G''$
implies the desired assertion.

	\qed(for the logical principle lemma.)

Clearly then  we get that , if  $c_i = G.C.D. (b_i , d_i ) $
and $\De$  is the
group  $\Ga ( m, b_1 / c_1 , b_2 / c_2 , \dots b_r / c_r ; d_1/ c_1
,d_2 / c_2 , \dots  d_r/
c_r )$, an iterated application of the logical principle yields that the
order of $\ga$ is the same in $\Ga$  and in $\De$.

{\bf  Step II.}

   Let $T : = T( m, b_1 , b_2 , \dots b_r ; d_1 ,d_2 , \dots d_r )$ be
the  quotient of the group
$\Ga ( m, b_1 , b_2 , \dots b_r , d_1 ,d_2 , \dots d_r ), $ by the
central cyclic subgroup
$C(\ga)$  generated by  $\ga$ : then by step I $T$ is isomorphic to 
the polygonal
group
$T( b_1 , b_2 , \dots b_r)$
with generators
$\de_1 ,\de_2 , \dots \de_r $, and relations $\de_1^{b_1} =\de_2^{b_2} =
   \dots  = \de_r ^{b_r} = \de_1 \cdot \de_2  \cdots \de_r = 1$.

In fact $ T( m, b_1 , b_2 , \dots b_r ; d_1 ,d_2 , \dots d_r )$ is a quotient
of the free product of cyclic groups of respective orders $b_i$ by the
relation that be trivial the product  $\beta_1 ^{d_1} \cdot
\beta_2 ^{d_2}  \cdots \beta_r ^{d_r}$.  But ,
since  $G.C.D. (b_i , d_i )= 1$ , each  $ \beta_i ^{d_i}  : = \de_i$
is a generator of the respective
cyclic group.

   {\bf  Steps III-V .}

We have thus a central extension

$$1 \ra C(\ga) \ra  \Ga ( m, b_1 , b_2 , \dots b_r ; d_1 ,d_2 , \dots d_r )
    \ra T(b_1 , \dots b_r )  \ra 1 ,$$
where $C(\ga)$ is the cyclic central subgroup generated by $\ga$ , and the
quotient $T: = T(b_1 , \dots b_r ) $ is the polygonal group defined above.

Our strategy will consist in  proving  that  either

III) the image of $\ga$  is nontrivial in $\Q$-homology ( i.e., in the
Abelianization of $\Ga$ tensored with $\Q$), whence a fortiori $\ga$
has infinite
order in $\Ga$, or

IV) $H^1( \Ga, \Q ) = 0$ : however then , in the nonexceptional cases,
$\Ga$
differs from $T$  because it has cohomological dimension 3 instead of 2, and
thus in any case $\ga$  has infinite order in $\Ga$.

V) treats then the exceptional cases using integral homology and  matrix
representations.

{\bf Step III.}

The above odd looking alternative is a consequence of the following

\begin{prop}\label{homology}
Let $\Ga$ be the above group $\Ga ( m, b_1 , b_2 , \dots b_r ; d_1
,d_2 , \dots d_r )$.
   Then then the image of $\ga$  in $H_1( \Ga, \Q )$ is a generator, and it
is nonzero if and only if $ m  \neq  \sum_i  ( d_i / b_i ).$
\end{prop}

\Proof. Let  $[\ga] , [\beta_i ]$ be the respective images of  $\ga ,
\beta_i$, inside
$H_1( \Ga, \Q )$. Then they generate it and there are only the relations
$$[\beta_i ] =
(1/ b_i ) [\ga]  , {\rm and } \  ( m - \Sigma_i  ( d_i / b_i ) ) [\ga]= 0. $$
Whence, $[\ga]$ generates
$H_1( \Ga, \Q )$ and $H_1( \Ga, \Q ) \neq 0 $ if and only if
   $ m  =  \Sigma_i  ( d_i / b_i ).$

\qed

{\bf Step IV.}

Assume then that $H_1( \Ga, \Q ) = 0$ , and observe that, because of
our plumbing
construction, $\Ga$ is the fundamental group of an orientable
3-manifold $M : = \Sigma$.
In particular, $H_1( M, \Q )=  H_1( \Ga, \Q )= 0 $,  and by
Poincar\'e Duality and
ordinary duality  $H^1( M, \Q )= H^2 ( M, \Q )= 0$, while  $H^3 ( M,
\Q )\cong \Q.$
   Let $N$ be the universal covering of $M$: then we have a spectral sequence
$ H^p ( \Ga, H^q ( N, Q ))$ converging to the graded module
associated to a suitable
filtration of  $H^{p+q} ( M, Q )$, for each ring $Q$ ($ Q = \Z$ or
$\Q $ in our
application).

Clearly, $H^1 ( N, \Q )=0$, hence $H^2 ( M, \Q )=0 $  implies $H^2 (
\Ga, \Q )=0.$

We can moreover apply ( cf. \cite{wei} 6.8.2.) the Lyndon-Hochshild-Serre
spectral sequence associated to the exact sequence

$1 \ra C(\ga) \ra  \Ga : = \Ga ( m, b_1 , b_2 , \dots b_r ; d_1 ,d_2
, \dots d_r )
    \ra T  \ra 1 ,$

whose $E_2$ term is $H^p ( T , H^q ( C(\ga), \Q ))$ and which
converges to a graded
quotient of $H^{p+q} ( \Ga, \Q ).$

   Now, if $\ga$  had finite order, then $H^i (
C(\ga), \Q )=0 $ for each $i \geq 1$, whence $H^i ( \Ga, \Q )= H^i ( T, \Q ) $
for each $i \geq 0$.

We get therefore an obvious contradiction in the case where $H^2 ( T, \Q ) \neq
0. $

   Observe that the polygonal group $T$ is a quotient of the group $\Pi$ with
generators $\beta_1 ,\beta_2 , \dots \beta_r $ , and with relation
$\beta_1 \cdot \beta_2  \cdots \beta_r  =1$. $\Pi$   is the
fundamental group of $\PP^1_{\C}$ minus $r$ points, and $T$ is the orbifold
fundamental group of the maximal Galois cover $C$ of $\PP^1$ branched in these
points with respective ramification multiplicities exactly
equal to $b_1 -1, b_2 -1, \dots
b_r -1.$

If $T$ is infinite, then $C$ is not compact, otherwise  $ C
\cong\PP^1$ , by the Riemann
mapping theorem.   Whence if $T$ is infinite, $H^2 (\PP^1 , \Q )
\cong \Q \cong H^2 ( T, \Q
) $ and we have found the required contradiction.

Otherwise, $T$ is finite,
and $ C  \ra \PP^1 $ has a finite degree $d$ .  As well known, by the
formula of
Hurwitz, then  $2 - 2/d = \Sigma_i ( 1- 1 / b_i )$ which implies that
$r \leq 3$, and
since $ r \geq 3$ we get $r=3 $ and  $ \Sigma_i ( 1- 1 / b_i ) > 1$,
an inequality which leads
us to the exceptional cases for $(b_1, b_2, b_3)$, corresponding to
the Platonic solids and to the Klein groups

Va ) $( 2,2,n)$ , $ n \geq 2$ ,( $d= 2n$) , $(d_1 , d_2 , d_3 ) =(1,1,t) $

Vb )  $( 2,3, n)$, $ 3 \leq n \leq 5$ ( $d = 12, 24, 60$),  $(d_1 ,
d_2 , d_3 ) =(1,c,t) $

(here $ c= 1,2$ and $ 1 \leq t \leq n-1$).

{\bf Step Va.}

Assume we are in the exceptional case a): in this case we shall explicitly
prove that the group $\Ga$ is finite, find a faithful matrix representation,
and find that the period of $\ga$ equals exactly $2p$, where $ p : = 
(m-1)n -t$.
Thus, the order of $\ga$ is always $ \geq 2$.

   In fact, we can
change the presentation of the group, eliminating $ \ga  = \beta_3
^{b_3} = \beta_3^n$   and
obtaining the relation  $\beta_3 ^{mn - t} =  \beta_1  \cdot \beta_2$ .

Then, $ \beta_1  \cdot \beta_2= \beta_3 ^{mn - t} =
\beta_1^2  \cdot \beta_3^p$, whence  $\beta_2 = \beta_1  \cdot \beta_3^p$.

Setting for simplicity  $a : = \beta_1 , b : =
\beta_3$ , we get the presentation
$$\Ga = \langle a,b | a^2  = b^n= a \cdot b^p \cdot a \cdot b^p \rangle.$$

   Since
$a^2  =   a \cdot b^p \cdot a \cdot b^p$  , we get  $ b^{-p}= a b^p
a^{ -1} $, whence
     $b^{-pn}= a b^{pn} a^{ -1 }$ and since $a$ commutes with  $b^n=
a^2$  , finally that
$b^{-pn}=  b^{pn} $  , i.e., $ b^{2pn}= 1= a^{4 p}$ .

It follows that the order of the
group $\Ga$ is at most  $4 pn$, and that equality holds if the period of $a$ is
exactly equal to $4p$.

In order to show that the period of $a$ is
exactly equal to $4p$ we use  the following representation
$\rho : \Ga \ra   GL(2,\C),$
such that
$$\rho (a)=\left(
\begin{array}{cc}
0 & \zeta _{4p}\\
\zeta _{4p}& 0
\end{array}
\right).
$$
$$\rho (b)=\left(
\begin{array}{cc}
\zeta _{2np} & 0 \\
0 &  u \ \zeta _{2np}^{-1}
\end{array}
\right).
$$
where $\zeta _h$  is  $: = exp (
2 \pi i / h) $, and $u$ is a  $p$-th root of $1$ such that $u^n=\zeta_p$
(recall that, since we assumed $G.C.D. \ (n,t)=1$, also $G.C.D. \ (p,n)=1$).

One can indeed
verify that $\rho (a^2 )= \rho (b^n) =
\rho ( (a \cdot b^p)^2) = \zeta _{2p}   \cdot  Id$, as claimed.

						\qed

{\bf Step Vb.}

Assume that we are in the exceptional case b).

In this case, we shall first try
  to  show that the image of $\ga$ in the abelianization $G$ of
$\Ga$ is nontrivial.

Eliminating $\ga$ we get  $\beta_1  = \beta_3 ^{mn - t }
\beta_2^{-c}  $, thus $\Ga$ is generated by  $a: =  \beta_2  , b:= \beta_3 $,
with relations
$$a^3
= b^n=  b^{p+n} a ^{-c}  b^{p+n} a ^{-c }, $$
where $p: = n(m-1) -t $, as above.

Letting $A,
B$ , be the respective images of $a,b$, in the abelianization of $\Ga$, we
obtain:
$$3 A - n B = 0 , 2c A =( 2 p + n )B.$$
   Since $3 - 2c =  \pm 1$ ( according to the respective cases  $c=1 , 
c=2$), we
get the relation  $\pm A + 2p  B = 0$,  thus $G$ is cyclic with generator $B$.

Moreover, the relation $ n B = 3 A = -( \pm 6 p  B )$ shows that $B$ 
has period
$f : = n  \pm 6p .$

    Now, if  $m \geq  2$, then $ p > 0$, thus if $c=1$
then $ f > n$, whence  $n B \neq  0$, as we wanted to show.

   If instead $m \geq  2, c= 2 $, the absolute value of the period equals
$ 6 p - n = n [6
(m-1) -1] = 6t $, which is clearly $ > n$  as soon as $ m \geq 3$.

If instead $ m =2$, the absolute value of the period is $ > n$ iff
$ 4 n > 6 t $, which holds unless $ \frac23 \  n \leq t \leq (n-1)$,
i.e., unless $ t = n-1$.

But in this case one has $ f = 5 n - 6 (n-1) = 6 -n$, thus
$ n B = 0$ since $6-n$ divides $n$.

Similarly, if $ m=1, p = -t$, we have $ f = \pm  6 t - n$, and $n B \neq 0$ if
$c= 2$, whereas if  $c=1$ we can reach this conclusion only if $n$ is not a
multiple of  $6t-n $.

This condition then holds unless $ t = 1$, and $ n = 3,4,5$.

We are left then with two cases to consider, the first where
$c= 2$, the second where $c=1$. For the latter case, we 
use directly
a result which goes back essentially to Felix Klein (\cite{klein}), 
and is clearly
stated by Milnor in \cite{milnor}:

Given a triangle group $ T:= T(1, b_1, b_2, b_3; 1,1,1)$ which is 
elliptic, i.e.,
such that $ \Sigma_i \frac{1}{b_i} > 1$,
then its inverse image $\hat{T}$ in $ SU(2, \C)$ has the presentation
$$ \hat{T} = \langle \ga, \beta_1 ,\beta_2 ,  \beta_3 | \ga = \beta_1 ^{b_1}=
\beta_2 ^{b_2} =   \beta_3 ^{b_3} = \beta_1  \cdot
\beta_2   \cdot \beta_3 \rangle .$$

It follows that $\hat{T}$ is isomorphic to our group $\Ga$, thus we have a
nontrivial central extension of $T$ by the central element $\ga$ of
of order two.

In the former case, we have the following presentation for $\Ga$
$$ \Ga = \langle \ga, \de_1 ,\de_2 ,  \de_3 | \ga = \de_1 ^{2}=
\de_2 ^{3} =   \de_3 ^{n} , \ga^2 = \de_1  \cdot
\de_2   \cdot \de_3^{n-1} \rangle .$$

Again here we use the extended triangle group $\hat{T}$,
setting
$$\de_1 : = \beta_1, \de_2  : = \beta_2 , \de_3 := \beta_3^{-1}. $$
Then we see that we get a homomorphic image of $\Ga$, where $\ga$ maps
onto an element  of order 2 (that we still denote by $\ga$).

We are finished with Vb).

				\qed

\bigskip

\section{Proofs of the main theorems}

\Proof {\em of Theorem A}
By remark \ref{simplified} we may replace $\Ga$ by its homomorphic 
image given by
the simplified group. I.e., we may assume $ g_i = 1$ or $ = 0$.

If $g_i \geq 1$, by remark \ref{arbitrary}, we may again take a homomorphic
image of $\Ga$ corresponding to changing $ g_i$ to $0$, and to changing $m_i$
making it arbitrarily high (i.e., making the self-intersection 
extremely negative).

Thus we may assume that we have a tree of rational curves, where
$ - m_i \leq -2, \ \forall i$.

If the tree is linear, the statement follows by lemma \ref{ltree}.

If we have a comb of rational curves, and $\ga_i$ corresponds to the rim of the
comb,
then the nontriviality of $\ga_i$ follows by theorem \ref{rim} and by the
subsequent Steps III, IV, V; else, it follows by proposition
\ref{combred}.

The remaining cases are taken care of, again by proposition \ref{combred}.

\qed

\Proof {\em of Theorem B}

Observe that if ii-1) holds, and $g_i = 0$, then if $D'_i$ is a curve
we have $ K_{S'} \cdot D'_i \geq 0$, hence also $ K_{S} \cdot D_i \geq 0$.

Thus we see that all the curves $D_i$ with $g_i = 0$ have self-intersection
$ D_i ^2 = - m_i \leq -1$, therefore assumption ii-1) implies assumption
ii-2) and we proceed with assumption ii-2), without forgetting the other
assumption of minimality in the GNC category. This implies that if
$g_i = 0$ and $ D_i ^2 = - 1 $, then $D_i$ meets at least three other 
components.

We can then use exactly the same strategy used for theorem A, since
the case of a linear tree follows automatically, and curves with
self-intersection $-1$ occur only as rims, and in this case the possibility
$ m=1$ is contemplated in theorem \ref{rim} and in the
subsequent Steps III, IV, V.

\qed

\Proof {\em of Theorem C}

We follow again the strategy of proof of theorem A.

If we have a linear tree, and there is a curve of positive genus, then we may
conclude that each $\ga_i$ has infinite order by remark \ref{infty}.

If we have a comb, then we know by theorem \ref{rim} that the generator $\ga$
corresponding to the rim has infinite order, if we are not in the 
exceptional cases
Va), Vb). Let moreover $\ga_i$ belong, say, to the string $ D(1)$.

Then we have shown in  5A (cf. lemma \ref{ltree}) that $ \ga = \ga_1^{a_n}$,
and $ \ga_i = \ga_1^{a_i}$, where $ 1 \leq a_i \leq a_n$.

Hence, also $ \ga_1$ and $ \ga_i$ have infinite order in 
the nonexceptional cases.

Similarly we are done if we have a comb and there is a curve $D_i$ of positive
genus, since we may then reduce to the case where all the genera are $0$, but
$ m_i$ is arbitrary, hence we are not in the exceptional cases.

So our statement is proven for elementary infinite pieces, and the rest
  follows easily by induction, since we may apply lemma \ref{G1G2} and corollary
\ref{3ormore}.

\qed

\begin{footnotesize}
\noindent
{\bf Note .}
When I presented these results at the AMS Meeting in NY,
november 3-5 2000, Walter Neumann mentioned that our presentation
of the local  fundamental group of neighbourhoods of divisors in
complex surfaces
is  similar to the method of \cite{neu81} of solid tori
   decompositions for 3-manifolds (in turn
   based on the methods  earlier
introduced by Waldhausen (\cite{wald}, \cite{wald2}),
who studied the problem whether  such manifolds
would be determined by their fundamental group.

We would also like to mention that Wagreich (\cite{wag}) and Karras
(\cite{kar}) determined the cases where $D$ comes from a singularity
and the group $\Ga$ is solvable.

\end{footnotesize}

\bigskip
\noindent
{\bf Acknowledgements.}
I would like to  express my indebtness
to Guido Zappa, for a precious suggestion and for the
letters exchange I had with him, and also to Nick Shepherd-Barron,
for another precious suggestion that came out during our train
conversation, in june 98 near Oberwolfach. Thanks to Igor Dolgachev
for kindly pointing out the  relevance of reference \cite{milnor}.

\vfill

\noindent
{\bf Author's address:}

\bigskip

\noindent
Prof. Fabrizio Catanese\\
Lehrstuhl Mathematik VIII\\
Mathematisches Institut\\
NW II, Zi. 730,\\
Universit\"at Bayreuth\\
   D- 95448 BAYREUTH, GERMANY

e-mail: Fabrizio.Catanese@uni-bayreuth.de

\end{document}